\documentclass[a4,12pt,leqno]{article}
\usepackage{amsmath}
\usepackage{graphics}
\usepackage{latexsym}
\usepackage{amsfonts}
\numberwithin{equation}{section} \rightmargin 1.5cm \leftmargin
-1.5cm \textheight 23cm \textwidth 16cm \topmargin -1cm
\title{
 Minimax Impulse Control Problems\\ in Finite Horizon
}
\author{Brahim EL ASRI \thanks
{%
Universit\'e Cadi Ayyad, D\'ept. de Math\'ematiques,  FSTG, B.P.
549, Marrakech, 40.000, Maroc. e-mails: b.elasri@uca.ma  }}
\begin{document}
\date{}
\maketitle
\newtheorem{theo}{Theorem}
\newtheorem{problem}{Problem}
\newtheorem{pro}{Proposition}
\newtheorem{cor}{Corollary}
\newtheorem{axiom}{Definition}
\newtheorem{rem}{Remark}
\newtheorem{lem}{Lemma}
\newcommand{\brm}{\begin{rem}}
\newcommand{\erm}{\end{rem}}
\newcommand{\beth}{\begin{theo}}
\newcommand{\eeth}{\end{theo}}
\newcommand{\bl}{\begin{lem}}
\newcommand{\el}{\end{lem}}
\newcommand{\bp}{\begin{pro}}
\newcommand{\ep}{\end{pro}}
\newcommand{\bcor}{\begin{cor}}
\newcommand{\ecor}{\end{cor}}
\newcommand{\be}{\begin{equation}}
\newcommand{\ee}{\end{equation}}
\newcommand{\beq}{\begin{eqnarray*}}
\newcommand{\eeq}{\end{eqnarray*}}
\newcommand{\beqa}{\begin{eqnarray}}
\newcommand{\eeqa}{\end{eqnarray}}
\newcommand{\dg}{\displaystyle \delta}
\newcommand{\cm}{\cal M}
\newcommand{\cF}{{\cal F}}
\newcommand{\cR}{{\cal R}}
\newcommand{\bF}{{\bf F}}
\newcommand{\tg}{\displaystyle \theta}
\newcommand{\w}{\displaystyle \omega}
\newcommand{\W}{\displaystyle \Omega}
\newcommand{\vp}{\displaystyle \varphi}
\newcommand{\ig}[2]{\displaystyle \int_{#1}^{#2}}
\newcommand{\integ}[2]{\displaystyle \int_{#1}^{#2}}
\newcommand{\produit}[2]{\displaystyle \prod_{#1}^{#2}}
\newcommand{\somme}[2]{\displaystyle \sum_{#1}^{#2}}
\newlength{\inter}
\setlength{\inter}{\baselineskip} \setlength{\baselineskip}{7mm}
\newcommand{\no}{\noindent}
\newcommand{\rw}{\rightarrow}
\def \ind{1\!\!1}
\def \R{I\!\!R}
\def \N{I\!\!N}
\def \cadlag {{c\`adl\`ag}~}
\def \esssup {\mbox{ess sup}}
\begin{abstract}
We consider the problem of impulse control minimax in finite
horizon, when cost functions $(C(t,x,\xi)>0)$. We show existence of
value function of the problem. Moreover, the value function is
characterized as the unique viscosity solution of an Isaacs
quasi-variational inequality. This problem is in relation with an
application in mathematical finance.
%In this paper we show existence and uniqueness of a solution for a
%system of $m$ variational partial differential inequalities with
%inter-connected obstacles. This system is the deterministic version
%of the Verification Theorem of the Markovian optimal $m$-states
%switching problem. The switching cost functions are arbitrary. This
%problem is in relation with the valuation of firms in a financial
%market.
\end{abstract}

%\no{\bf AMS Classification subjects}: 60G40 ; 62P20 ; 91B99 ; 91B28
%; 35B37 ; 49L25.
\medskip

\no {$\bf Keywords$}: Impulse control; Robust control; Differential
games; Quasi-variational inequality; Viscosity solution

\section {Introduction}
%We consider an optimal impulse control problem with finite horizon.
%In the game, player-$\xi$ would like to minimize the pay-off by
%choosing suitable impulse control $\xi(.)$, whereas player-$\tau$
%wants to maximize the pay-off by choosing a proper control. In
%mathematical finance, one may consider the option pricing problem of
%references \cite{[BE], [BET]}.  If the piecewise linear transaction
%costs are replaced by a more realistic piecewise affine cost, i.e. a
%fixed cost is charged for any transaction in addition to a variable
%part, then the problem at hand is exactly that considered here.\\
In this paper we study an optimal impulse control problem with
finite horizon.\\ \indent Optimal impulse control problems appear in
many practical situations. In the game, player-$\xi$ would like to
minimize the pay-off by choosing suitable impulse control $\xi(.)$,
whereas player-$\tau$ wants to maximize the pay-off by choosing a
proper control. In mathematical finance, one may consider the option
pricing problem of references \cite{[BE], [BET]}.  If the piecewise
linear transaction costs are replaced by a more realistic piecewise
affine cost, i.e. a fixed cost is charged for any transaction in
addition to a variable part, then the problem at hand is exactly
that considered here. We refer the reader to \cite{[BL]} (and the
references cited therein) for extensive discussions. For
deterministic autonomous systems with infinite horizon, optimal
impulse control problems were studied in \cite{[B]}, and optimal
control problems with continuous, switching, and impulse controls
were studied by the author \cite{[Y11]} (see also \cite{[Y1]}).
Differential games with switching strategies in finite and infinite
duration were also studied \cite{[Y2], [Y3]}.  J. Yong, in
\cite{[Y]}, also studies differential games where one person uses an
impulse control and
other uses continuous controls.\\
\indent The study of optimal control problems with continuous
controls, gives rise to Hamilton-Jacobi-Bellman equations which are
satisfied by the value function corresponding to the problem, if it
is smooth.  It is known that the value function of these problems,
whenever smooth satisfy different variational and quasivariational
inequalities. But most of the time these value functions are only
continuous and not sufficiently smooth. The notion of viscosity
solutions, a kind of generalized solutions, introduced by Crandall
and Lions \cite{[CIL]} is extremely well suited for these problems.
The value function satisfies the corresponding equations or
inequalities in the viscosity sense. These control problems are
studied in the viscosity solution set up, for example, in \cite{[B],
[CE]}.  In all these works the existence and uniqueness results are
obtained assuming that the dynamics and cost functionals are bounded
and uniformly continuous and hence the value functions are in the
bounded uniformly continuous function class.

 In the finite horizon framework El Farouq et al
\cite{[EBB]} extended the work of Yong \cite{[Y]} but allowing
general jumps. In this works the existence of the value functions of
optimal impulse control problem and uniqueness of viscosity solution
are obtained assuming that the dynamics and costs functionals are
bounded and the impulse cost function should not depend on $x$.
Recently El Asri \cite{[E1]} have considered this
 impulse control problem when the dynamics unbounded and
costs functionals are bounded from below and the impulse cost
function depends on $x$.\\
\indent The purpose of this work is to fill in this gap by providing
a solution to the optimal impulse control problem using dynamic
programming principle tools and partial differential equation
approach.\\
 \indent We prove existence of the value function of the
problem when costs functionals  $C> 0$. We show that the value
function of the problem is associated  of deterministic functions
$v$ which is the unique solution of the following system of PDIs
(Isaacs equation):

\be
\begin{array}{l}\label{Issacs1}
\displaystyle \max\left\{\min\limits_{\tau\in{\cal
K}}\left[-\frac{\partial v}{\partial t}- \frac{\partial v}{\partial
x}f(t,x,\tau)-\psi(t,x,\tau)\right],\right.\\\qquad\qquad \qquad
\left.v(t,x)-N[v](t,x)\right\}= 0,
\end{array}
\ee where $\mathcal{K}$ compact. It turns out that this Isaacs
equation is the deterministic version of optimal impulse control
problem in finite horizon.\\
\indent This paper is organized as follows:\\
\indent In Section 2, we formulate the problem, we give the related
definitions and we prove the value function is bounded from below
with linear growth.  In Section 3, we give some properties of the
value function, especially the dynamic programming principle. Then
we introduce the approximating scheme which enables us to construct
a solution for value function of the optimal impulse control problem
which, in combination with the dynamic programming principle, play a
crucial role in the proof of the existence of the value function.
Section 4 is devoted to the connection between the optimal impulse
control problem and quasi-variational inequality. In Section 5, we
show that the solution of  QVIs is unique in the subclass of bounded
from below continuous functions which satisfy a linear growth
condition. $\Box$
\section{Formulation of the problem and preliminary results}

\subsection{Setting of the problem}
Let a two-players differential game system be defined by the
solution of following dynamical equations
\begin{equation}\label{dynamic}\left\{
\begin{array}{l}
\dot{y}(t)=f(t,y(t),\tau(t))\\
y(t_0)=x\in \R^m,\\
y(t_k^+)=y(t_k^-)+g(t_k,y(t^-_k),\xi_k),\,\,t_k\geq t_0,\,\xi_k\neq
0,%\\
%y(t^-_k)=y(t_k),
\end{array}
\right.
\end{equation}
where $y(t)$ is the state of the system, with values in $\R^m$, at
time $t$, $x$ is the initial state. The time variable $t$ belongs to
$[t_0,T]$ where $0\leq t_0<T$, and
$y(t_k^{\pm})=\lim\limits_{t\rightarrow t_k^{\pm}}y(t)$. We assume
that $y$ is left continuous at the times
$t_k$: $y(t_k^-)=y(t_k)$, $k\geq 1$. \\
The system is driven by two controls, a continuous control
$\tau(t)\in \mathcal{K}\in \R^m$, where $\mathcal{K}$ is compact
set, and an impulsive control defined by a double sequence $t_1, . .
. , t_k, . . . , \xi_1, . . . , \xi_k, . . .,\\k \in \N^* =\N
\backslash \{0\}$, where $t_k$ are the strategy, $t_k\leq t_{k+1}$
and $\xi_k\in\R^m$ the control at time $t_k$ of the jumps in
$y(t_k)$. Let $\mathcal{S}:=((t_k)_{k\geq 1},(\xi_k)_{k\geq 1})$ the
set of these strategies denoted by
$\cal D$.\\
For any initial condition $(t_0, x)$, controls $\tau(\cdot)$ and
$\cal S$ generate a trajectory $y(\cdot)$ of this system. The
pay-off is given by the following:

\begin{equation}\label{pay-off} J(t_0,x,\mathcal{S},\tau(.))=\integ{t_0}{T}\psi(s,y(s),\tau(s))ds
+\sum_{k\geq 1} C(t_k,y(t_k),\xi_k)\ind_{[t_{k}\leq T]}+G(y(T)),
\end{equation}
where if $t_k=T$ for some $k$ then we take $G(y(T))=G(y(T^+))$. The
term $C(t_k,y(t_k),\xi_k)$ is called the impulse cost. It is the
cost when player-$\xi$ makes an impulse $\xi_k$ at time $t_k$. In
the game, player-$\xi$ would like to minimize the pay-off
(\ref{pay-off}) by choosing suitable impulse control $\xi(.)$,
whereas player-$\tau$ wants to maximize the pay-off (\ref{pay-off})
by choosing a proper control $$\tau(.)\in \Omega=\{\mbox{measurable
functions}\quad [t_0,T]\rightarrow \mathcal{K}\}.$$ We shall
sometimes write $\tau\in \Omega$ instead of $\tau(.)\in \Omega$.\\

We now define the admissible strategies $\varphi$ for the minimizing
impulse control $\cal D$, as non-anticipative strategies. We shall
let ${\cal D}_a$ be the set of all such non-anticipative strategies.
\begin{axiom}
A map $\varphi: \Omega\rightarrow \mathcal{S}$ is called a
non-anticipative strategy if for any two controls $\tau_1(.)$ and
$\tau_2(.)$, and any $t\in[t_0,T ]$, the condition on their
restrictions to $[t_0, t[: \tau_1|_{[t_0,t [} = \tau_2|_{[t_0,t [} $
implies $\varphi(\tau_1)|_{[t_0,t ]} = \varphi(\tau_2)|_{[t_0,t ]}$.
\end{axiom}

In the next, we define the value function of the problem
$v:[0,T]\times \R^m\rightarrow \R$ as
$$v(t_0,x)=\inf\limits_{\varphi\in {\cal D}_a}\sup\limits_{\tau(.)\in \Omega}J(t_0,x,\varphi(\tau(.)),\tau(.))$$

\subsection{Assumptions} Throughout
this paper $T$ (resp. $m$) is a fixed real (resp. integer) positive
constant. Let us now consider the followings:
\medskip

\indent $(1)$ $f:[0,T]\times\R^m\times \mathcal{K}\rightarrow
\R^{m}$ and $g:[0,T]\times\R^m \times \R^m\rightarrow \R^{m}$ are
two continuous functions for which there exists a constant $C\geq 0$
such that for any $t\in [0, T]$, $\tau \in \mathcal{K}$ and $\xi, x,
x'\in \R^m$\be \begin{array}{ll}\label{regbs1}|f(t,x,\tau)|+
|g(t,x,\xi)|\leq
C(1+|x|) \quad \mbox{ and } \quad \\
|g(t,x,\xi)-g(t,x',\xi)|+|f(t,x,\tau)-f(t,x',\tau)|\leq
C|x-x'|\end{array}\ee

$(2)$ $C:[0,T]\times\R^m \times \R^m\rightarrow \R^{}$, is
continuous with respect to $t$ and $\xi$ uniformly in $x$.\\
For any $(t,x,\xi) \in [0,T]\times\R^m \times \R^m$, $C(t,x,\xi)$
are satisfying \be \label{CC}C(t,x,\xi) > 0. \ee
 %with
%linear growth $$|C(t,x,\xi)|\leq C(1+|x|),\, \, \forall (t,x,\xi)\in
%[0,T]\times \R^m\times  \R^m.$$
%Moreover we assume that \be \label{CCC}C(t,x,\xi)\rightarrow +\infty
%\quad \mbox{as} \quad \xi\rightarrow +\infty. \ee

$(3)$ $\psi:[0,T]\times\R^m\times \mathcal{K}\rightarrow \R$ is
continuous  with respect to $t$ and $\tau$ uniformly in $x$ with
linear growth, \be \label{polycond} |\psi(t,x,\tau)|\leq C(1+|x|),\,
\, \forall (t,x,\tau)\in [0,T]\times \R^m\times  \mathcal{K},\ee and
is bounded from below.

$(4)$ $G:\R^m\rightarrow \R$ is uniformly continuous  with linear
growth \be \label{polycond1}|G(x)|\leq C(1+|x|),\, \, \forall x\in
 \R^m, \ee
 and
is bounded from below.
%\\
%with respect to $t$ and
%$\tau$ uniformly in $y$, there exist some positive constants $C$ and
%$\gamma$ such that: \be \label{polycond}|G(x)|+
%|\psi(t,x,\tau)|+|C(t,x,\xi)|\leq C(1+|x|^\gamma),\, \, \forall
%(t,x,\tau,\xi)\in [0,T]\times \R^m\times  \mathcal{K} \times \R^m.
%\ee Moreover we assume that there exists a constant $\alpha >0$ such
%that for any $(t,x,\xi)\in [0,T]\times \R^m\times \R^m$, \be
%C(t,x,\xi)\geq \alpha. \ee %This condition means that
%switching from one mode to another one is not free and costs at
%least  $\alpha>0$.
%
%$(5)$ $\psi$ and $G$ are bounded.

These properties of $f$ and $g$ imply in particular that $y(t)_{0\le
t\leq T}$ solution of the standard DE (\ref{dynamic}) exists and is
unique, for any $t\in [0, T]$ and $x\in \R^m$.
\medskip

\subsection{Preliminary results}
We want to investigate the problem of minimizing
$\sup\limits_{\tau\in \Omega}J$ through the impulse control. We mean
to allow closed loop strategies for the minimizing control. We
remark that, being only interested in the $\inf \sup$ problem, and
not a possible saddle point.
\begin{theo}
Under the standing assumptions (Sect. 2.2) the value function $v$ is
bounded from below with linear growth.
\end{theo}
$Proof$: Consider the particular strategy in ${\cal D}_a$ is the one
where we have no impulse time. In this case $\sum_{k\geq 1}
C(t_k,y(t_k),\xi_k)\ind_{[t_{k}\leq T]}=0.$

$$ v(t,x)\leq \sup\limits_{\tau \in \Omega}[\integ{t}{T}\psi(s,y(s),\tau(s))ds +G(y(T))].
$$

Since $\psi$ and $G$ are linear growth, then
$$ v(t,x)\leq \integ{t}{T}C(1+|y(s)|)ds +C(1+|y(T)|).
$$
Now by using standard estimates from ODE, Gronwall's Lemma and the
strategy where we have no impulse time, we can show that
$$|y(t)|\leq C(1+|x|),$$
where $C$ is constant of $T$. Hence using this estimate we get

$$ v(t,x)\leq C(1+|x|).$$
 On the other hand, since the cost $C(t_k,y(t_k),\xi_k)$
are non negative functions and since $\psi$ and $G$ are bounded from
below,
then  $v$ is bounded from below.$\Box$\\

 We also state the following definition:
\begin{axiom}
For any function $v:[t_0,T]\times \R^m\rightarrow \R$, let the
operator $N$ be given by
$$N[v](t,x)=\inf\limits_{\xi\in E}[v(t,x+g(t,x,\xi))+C(t,x,\xi)].$$
\end{axiom}
%\begin{rem}
%If the function $v$ is continuous, so is the function $N[v]$.
%\end{rem}

\section{The value function}
\subsection{Dynamic programming principle}
The dynamic programming principle is a well-known property in
 optimal  impulse control. In our optimal control problem, it is
formulated as follows:
\begin{theo}(\cite{[EBB]}, Proposition 3.1)\label{principle-dyn}
The value function $v(.,.)$ satisfies the following optimality
principle:\\
for all $t\leq t'\, \in [t_0,T[$ and $x\in \R^m$,
$$v(t,x)=\inf\limits_{\varphi\in {\cal D}_a}\sup\limits_{\tau\in \Omega}[\integ{t}{t'}\psi(s,y(s),\tau(s))ds
+\sum_{k\geq 1,\,t_k<t'} C(t_k,y(t_k),\xi_k)\ind_{[t_{k}\leq
T]}+\ind_{[t' \leq T]}v(t',y(t'))],$$

and
$$v(t,x)=\inf\limits_{\varphi\in {\cal D}_a}\sup\limits_{\tau\in \Omega}[\integ{t}{t_n}\psi(s,y(s),\tau(s))ds
+\sum_{1\leq k<n} C(t_k,y(t_k),\xi_k)\ind_{[t_{k}\leq T]}+\ind_{[t_n
\leq T]}v(t_n,y(t_n))],$$ where $((t_n)_{n\geq 1},(\xi_n)_{n\geq
1})$ be an admissible control.
\end{theo}

\begin{pro}
The value function v(.,.) has the following property: \\
for all $t\in [t_0,T]$ and $x\in \R^m$,
$$v(t,x)\leq N[v](t,x).$$
\end{pro}
 $Proof:$
Assume first that for some $x$ and $t$:
$$v(t,x)> N[v](t,x).$$
Then we have for $t\leq t'$:\begin{eqnarray*} &&
\inf\limits_{\varphi\in {\cal D}_a}\sup\limits_{\tau\in
\Omega}[\integ{t}{t'}\psi(s,y(s),\tau(s))ds +\sum_{k\geq 1,\,t_k<t'}
C(t_k,y(t_k),\xi_k)\ind_{[t_{k}\leq T]}+\ind_{[t'\leq T]}v(t',y(t'))]\\
&& \quad >\inf\limits_{\xi\in E}[v(t,x+g(t,x,\xi))+C(t,x,\xi)].
\end{eqnarray*}
Among the admissible strategy $\varphi^\epsilon$'s  there are those
that place a jump at time t.
\begin{eqnarray*}
&& \sup\limits_{\tau\in \Omega}[\integ{t}{t'}\psi(s,y(s),\tau(s))ds
+\sum_{k\geq 1,\,t_k<t'}
C(t_k,y(t_k),\xi_k)\ind_{[t_{k}<T]}+\ind_{[t'<T]}v(t',y(t'))]\\
&& \quad > v(t,x+g(t,x,\xi))+C(t,x,\xi)-\epsilon.
\end{eqnarray*}
Now, pick $\tau_1$ such that
\begin{eqnarray*}
&& \integ{t}{t'}\psi(s,y(s),\tau_1(s))ds +\sum_{k\geq 1,\,t_k<t'}
C(t_k,y(t_k),\xi_k)\ind_{[t_{k}\leq T]}+\ind_{[t'\leq T]}v(t',y(t')) +\epsilon\\
&& \quad \geq \sup\limits_{\tau\in
\Omega}[\integ{t}{t'}\psi(s,y(s),\tau(s))ds +\sum_{k\geq 1,\,t_k<t'}
C(t_k,y(t_k),\xi_k)\ind_{[t_{k}\leq T]}+\ind_{[t'\leq
T]}v(t',y(t'))],
\end{eqnarray*}
which implies that:
\begin{eqnarray*}
&& \integ{t}{t'}\psi(s,y(s),\tau_1(s))ds +\sum_{k\geq 1,\,t_k<t'}
C(t_k,y(t_k),\xi_k)\ind_{[t_{k}\leq T]}+\ind_{[t'\leq T]}v(t',y(t')) +\epsilon\\
&& \quad > v(t,x+g(t,x,\xi))+C(t,x,\xi)-\epsilon.
\end{eqnarray*}
Choosing now $t'=t$, yields the relation
$$ \epsilon+v(t,x+g(t,x,\xi))> v(t,x+g(t,x,\xi))+C(t,x,\xi)-\epsilon.$$
By sending $\epsilon\rightarrow0$, we obtain $C(t,x,\xi)<0 $, which
is a contradiction.$\Box$
\subsection{Continuity of value function}
In this section we prove the continuity of the value function. The
main result of this section can be stated as follows.\\

We first present some preliminary results on $y(.)$. Consider
$\varphi \in {\cal D}_a$ and $(\tau(.),\varphi(\tau(.))$, composed
of jumps instants $t_1,t_2,...,t_n$ in the interval $[t,T]$, with
jumps $\xi_1,\xi_2,...,\xi_n,$ and let $y_1(.)$ and $y_2(.)$ be the
trajectories generated by ${\cal D}_a$, from $y_i(t)=x_i,\,i=1,2.$
\begin{lem}
There exists a constant $C$ such that for any $s\in [t,T]$,
$x_1,x_2\in \R^m$, and $k\in\{1,2...,n\}$
\begin{equation}\label{estimat2}
|y_1(s)-y_2(s)|\leq \exp(C(s-t))(1+C)^n|x_1-x_2|. \Box
\end{equation}
\end{lem}
$Proof:$ By the Lipschitz continuity of $f$ and Gronwall's Lemma, we
have

$$|y_1(s)-y_2(s)|\leq \exp(C(s-t))(1+C)|x_1-x_2|,\quad \forall s\in[t,t_1].$$
Next let us show for an impulse time
$$|y_1(t_k^+)-y_2(t_k^+)|\leq \exp(C(t_k-t))(1+C)^k|x_1-x_2|.$$
Looking more carefully at the first jump and using the Lipschitz
continuity of $g$, we have
\begin{equation}
\begin{array}{lll}
\label{viscder}
|y_1(t_1^+)-y_2(t_1^+)|&=|y_1(t_1^-)+g(t_1^-,y_1(t_1^-),\xi_1)-y_2(t_1^-)-g(t_1^-,y_2(t_1^-),\xi_1)|\\&\leq
(1+C)|y_1(t_1^-)-y_2(t_1^-)|
\\&\leq \exp(C(t_1-t))(1+C)|x_1-x_2|.
\end{array}
\end{equation}
The above assertion is obviously true for $k = 1$. Suppose now it
holds true at step $k$. Then, at step $k + 1$,
\begin{equation}
\begin{array}{lll}
\label{viscder} |y_1(t_{k+1}^+)-y_2(t_{k+1}^+)|&\leq
(1+C)|y_1(t_{k+1}^-)-y_2(t_{k+1}^-)|\\&\leq
(1+C)|y_1(t_k^+)-y_2(t_k^+)|\exp(C(t_{k+1}^- -t_k^+))
\\&\leq  \exp(C(t_{k+1}-t))(1+C)^{k+1}|x_1-x_2|.
\end{array}
\end{equation}
Finally $$|y_1(s)-y_2(s)|\leq
\exp(C(s-t))(1+C)^n|x_1-x_2|,\quad\forall s\in [0,T]. \Box$$

\medskip
We are now ready to give the main Theorem of this article. The value
function of the problem in which the controller chooses $n$ impulse
time is defined as

\begin{equation}\label{pay-off1} q^n(t_0,x)=\inf\limits_{\varphi\in {\cal D}_a}\sup\limits_{\tau(.)\in \Omega}[\integ{t_0}{T}\psi(s,y(s),\tau(s))ds
+\sum_{k=1}^{n} C(t_k,y(t_k),\xi_k)\ind_{[t_{k}\leq T]}+G(y(T))].
\end{equation}
We will denote the value of making no impulse by $q^0$, which we
define as

\begin{equation}\label{pay-off2} q^0(t_0,x)=\sup\limits_{\tau(.)\in \Omega}[\integ{t_0}{T}\psi(s,y(s),\tau(s))ds
+G(y(T))].
\end{equation}
Now, consider the following sequential optimal stopping problems:\\
for all $t\leq t'\, \in [t_0,T]$ and $x\in \R^m$,
\begin{equation}\label{approx}w^n(t,x)=\inf\limits_{(t',\xi_{t'})}\sup\limits_{\tau\in
\Omega}[\integ{t}{t'}\psi(s,y(s),\tau(s))ds +
C(t',y(t'),\xi_{t'})\ind_{[t'\leq
T]}+w^{n-1}(t',y(t')],\end{equation} where $w^0(t,x)=q^0(t,x)$.

\begin{pro}
(i) The sequence $(w^n)_{n\geq 0}$ converges decreasingly.\\
(ii) For $n\in \N$, we have that $q^n(t,x) = w^n(t,x)$, for all
$(t,x)\in [t_0,T]\times\R^m$.\\
(iii) For all $(t,x)\in [t_0,T]\times\R^m$, the decreasing sequence
 $(q^n(t,x))_{n\geq 0}$ converges:
$$\lim\limits_{n\rightarrow \infty}q^n(t,x)=v(t,x).$$
\end{pro}
$Proof:$ (i) We show by induction on $n\geq 0$, that  for each
$(t,x)\in [t_0,T]\times\R^m$
$$C(1+\mid x \mid) \geq w^n(t,x)\geq w^{n+1}(t,x).$$
For $n=0$ the property is obviously true, since it is enough to take
$t'= T$ in the definition of $w^1$ to obtain that $w^0\geq w^1$. On
the other hand taking into account that $\psi$ and $G$ are linear
growth, then

$$w^0(t,x) \leq C(1+\mid x \mid).$$
Suppose now that, for some $n$, we have
$$w^n(t,x) \geq w^{n+1}(t,x).$$
Replace $w^{n+1}$ by $w^n$ in the definition of $w^{n+2}$, to obtain
that $w^{n+1}(t,x)\geq w^{n+2}(t,x)$.\\
On the other hand, since the cost $C(t_k,y(t_k),\xi_k)$ are non
negative functions and since $\psi$ and $G$ are bounded from below,
then  by induction  on $n\geq 0$, $w^n$ is bounded from below.

\noindent (ii) Assume first that for some $(t,x)\in[0,T]\times \R^m$
$$w^n(t,x)> q^n(t,x)=\inf\limits_{\varphi\in {\cal D}_a}\sup\limits_{\tau(.)\in \Omega}[\integ{t}{T}\psi(s,y(s),\tau(s))ds
+\sum_{k=1}^{n} C(t_k,y(t_k),\xi_k)\ind_{[t_{k}\leq T]}+G(y(T))],$$
and let the difference be $2\epsilon.$ Choose an admissible strategy
$\psi^\epsilon$ that approximates the infimum in the r.h.s. up the
$\epsilon$. Then,
$$w^n(t,x)-\epsilon \geq \sup\limits_{\tau(.)\in \Omega}[\integ{t}{T}\psi(s,y(s),\tau(s))ds
+G(y(T))]  +\sum_{k=1}^{n} C(t_k,y(t_k),\xi_k)\ind_{[t_{k}\leq
T]}.$$ Since for any $(t,x,\xi)\in [0,T]\times \R^m\times \R^m$,
$C(t,x,\xi)\geq 0$ then we have:

$$w^n(t,x)-\epsilon \geq w^0(t,x),$$
a contradiction with $w^n$ is non increasing.\\
Assume the contrary that
$$\inf\limits_{(t',\xi_{t'})}\sup\limits_{\tau\in
\Omega}[\integ{t}{t'}\psi(s,y(s),\tau(s))ds +
C(t',y(t'),\xi_{t'})\ind_{[t'\leq T]}+w^{n-1}(t',y(t')]=w^n(t,x)<
q^n(t,x),$$ and let the difference be $2\epsilon.$ Choose a
$(t_1,\xi_1)$ that approximates $w^n(t,x)$ up the
$\frac{\epsilon}{n}$. Then,
$$\sup\limits_{\tau\in
\Omega}[\integ{t}{t_1}\psi(s,y(s),\tau(s))ds +
C(t_1,y(t_1),\xi_{t_1})\ind_{[t_1\leq
T]}+w^{n-1}(t_1,y(t_1)]-\frac{\epsilon}{n}\leq q^n(t,x)
-2\epsilon.$$ Now choose a $(t_2,\xi_2)$ that approximates
$w^{n-1}(t_1,y(t_1))$ up the $\frac{\epsilon}{n}$. Then we have:
$$\sup\limits_{\tau\in
\Omega}[\integ{t_1}{t_2}\psi(s,y(s),\tau(s))ds +
C(t_2,y(t_2),\xi_{t_2})\ind_{[t_2\leq
T]}+w^{n-2}(t_2,y(t_2)]-\frac{\epsilon}{n}\leq
w^{n-1}(t_1,y(t_1)).$$ It implies that
$$\begin{array}{ll}\sup\limits_{\tau\in
\Omega}[\integ{t}{t_1}\psi(s,y(s),\tau(s))ds +
C(t_1,y(t_1),\xi_{t_1})\ind_{[t_1\leq
T]}+\integ{t_1}{t_2}\psi(s,y(s),\tau(s))ds\\
\qquad\qquad\qquad\qquad +C(t_2,y(t_2),\xi_{t_2})\ind_{[t_2\leq
T]}+w^{n-2}(t_2,y(t_2))]-\frac{2 \epsilon}{n}\leq q^n(t,x)
-2\epsilon.\end{array}$$ Then
$$\begin{array}{ll}\sup\limits_{\tau\in
\Omega}[\integ{t}{t_2}\psi(s,y(s),\tau(s))ds +
C(t_1,y(t_1),\xi_{t_1})\ind_{[t_1\leq
T]}+C(t_2,y(t_2),\xi_{t_2})\ind_{[t_2\leq T]}\\
\qquad\qquad\qquad\qquad\qquad\qquad\qquad\qquad\qquad\qquad+w^{n-2}(t_2,y(t_2))]-\frac{2
\epsilon}{n}\leq q^n(t,x) -2\epsilon.\end{array}$$ Repeating this
procedure $n$ times, we obtain
$$\sup\limits_{\tau(.)\in \Omega}[\integ{t}{T}\psi(s,y(s),\tau(s))ds
+\sum_{k=1}^{n} C(t_k,y(t_k),\xi_k)\ind_{[t_{k}\leq
T]}+G(y(T))]-\frac{n\epsilon}{n}\leq q^n(t,x) -2\epsilon.$$ But
placing the infimum  in the l.h.s. of the last inequality leads to a
contradiction.\\

(iii) Since  $(q^n(t,x))_{n\geq 0}$ is a non-increasing sequence,
then \begin{equation}\label{ine1} \lim\limits_{n\rightarrow\infty}
(q^n(t,x))_{n\geq 0}\geq v(t,x), \quad (t,x)\in
[t_0,T]\times\R^m.\end{equation} Now we show that
$\lim\limits_{n\rightarrow \infty}q^n(t,x)\leq v(t,x) , \quad
(t,x)\in [t_0,T]\times\R^m.$\\
Recall the characterization of (\ref{pay-off}) that reads as:
$$v(t,x)=\inf\limits_{\varphi\in {\cal D}_a}\sup\limits_{\tau(.)\in \Omega}[\integ{t}{T}\psi(s,y(s),\tau(s))ds
+\sum_{k\geq 1} C(t_k,y(t_k),\xi_k)\ind_{[t_{k}\leq T]}+G(y(T))].$$
Fix an arbitrary $\epsilon>0$. Let $\varphi=((t_n)_{n\geq
1},(\xi_n)_{n\geq 1})$ belongs to $D_a$ such that
$$v(t,x)+\epsilon \geq \sup\limits_{\tau(.)\in \Omega}[\integ{t}{T}\psi(s,y(s),\tau(s))ds
+\sum_{k\geq 1} C(t_k,y(t_k),\xi_k)\ind_{[t_{k}\leq T]}+G(y(T))].$$
Since for any $(t,x,\xi)\in [0,T]\times \R^m\times \R^m$,
$C(t,x,\xi)\geq 0$ then we have:
$$v(t,x)+\epsilon \geq \sup\limits_{\tau(.)\in \Omega}[\integ{t}{T}\psi(s,y(s),\tau(s))ds
+\sum_{k=1}^{n} C(t_k,y(t_k),\xi_k)\ind_{[t_{k}\leq T]}+G(y(T))].$$
Then
$$v(t,x)+\epsilon \geq \inf\limits_{\varphi\in {\cal D}_a}\sup\limits_{\tau(.)\in \Omega}[\integ{t}{T}\psi(s,y(s),\tau(s))ds
+\sum_{k=1}^{n} C(t_k,y(t_k),\xi_k)\ind_{[t_{k}\leq
T]}+G(y(T))]=q^n(t,x),$$ And then $$v(t,x)+\epsilon \geq
\lim\limits_{n\rightarrow\infty}q^n(t,x).$$ As $\epsilon$ is
arbitrary then putting $\epsilon\rightarrow 0$ to obtain:
\begin{equation}\label{ine2}v(t,x)\geq \lim\limits_{n\rightarrow\infty}q^n(t,x). \Box\end{equation}
%$$\inf\limits_{\varphi\in {\cal D}_a}\sup\limits_{\tau(.)\in \Omega}[\integ{t}{T}\psi(s,y(s),\tau(s))ds
%+\sum_{k=1}^{n} C(t_k,y(t_k),\xi_k)\ind_{[t_{k}\leq T]}+G(y(T))]\leq
%q^n(t,x) -\epsilon.$$

 \beth The
value function $v:[0,T]\times \R^m\rightarrow \R$ is continuous in
$t$ and $x$.\eeth $Proof.$ The proof is divided in three steps. Let
us consider $\epsilon >0$ and $(t',x')\in B((t,x),\epsilon)$ %and let
%us consider the following set of strategies:$$ \tilde
%D_a:=\left\{(\delta,\xi)=((t_n)_{n\geq 1}, (\xi_n)_{n\geq 1}) \in
%{\cal D}_a \mbox{ such that } \forall n\geq 1, \ind_{[\tau_n\leq T]}
%\leq \frac{C(1+(\epsilon +|x|))}{n}\right\}.$$ From Proposition 1,
%the impulse control of
%non-anticipative strategy optimal $(\delta,\xi)$ belongs to $\tilde D_a$.\\
\textbf{Step 1. } First let us show that $q^n$ is upper
semi-continuous. Recall the characterization of $q^n$ that reads as
$$\begin{array}{ll}q^n(t,x)=\inf\limits_{\varphi\in  D_a}\sup\limits_{\tau\in \Omega}[\integ{t}{T}\psi(s,y(s),\tau(s))ds
+\sum_{k=1}^n C(t_k,y(t_k),\xi_k)\ind_{[t_{k}\leq T]}+
G(y(T))],\end{array}$$
$$\begin{array}{ll}q^n(t',x')=\inf\limits_{\varphi\in D_a}\sup\limits_{\tau\in \Omega}[\integ{t'}{T}\psi(s,y'(s),\tau(s))ds
+\sum_{k=1}^n C(t_k,y'(t_k),\xi_k)\ind_{[t_{k}\leq T]}
+G(y(T))].\end{array}$$ Fix an arbitrary $\epsilon^1>0$. Let
$\varphi=((t_n)_{n\geq 1},(\xi_n)_{n\geq 1})$ belongs to $ D_a$ such
that
\begin{eqnarray*}
&&\sup\limits_{\tau\in \Omega}[\integ{t}{T}\psi(s,y(s),\tau(s))ds
+\sum_{k=1}^n C(t_k,y(t_k),\xi_k)\ind_{[t_{k}\leq
T]}+G(y(T))] \\
&& \quad \leq \inf\limits_{\varphi\in D_a}\sup\limits_{\tau\in
\Omega}[\integ{t}{T}\psi(s,y(s),\tau(s))ds +\sum_{k=1}^n
C(t_k,y(t_k),\xi_k)\ind_{[t_{k}\leq T]}+G(y(T))]+\epsilon^1\\ &&
\quad =q^n(t,x)+\epsilon^1.
\end{eqnarray*}
Also,
$$\begin{array}{ll}q^n(t',x')\leq \sup\limits_{\tau\in
\Omega}[\integ{t'}{T}\psi(s,y'(s),\tau(s))\ind_{[s\geq t']}ds
+\sum_{k=1}^n C(t_k,y'(t_k),\xi_k)\ind_{[t_{k}\leq
T]}+G(y(T))].\end{array}$$ Now pick $\tau_1$ such that

\begin{eqnarray*}
&& \sup\limits_{\tau\in
\Omega}[\integ{t'}{T}\psi(s,y'(s),\tau(s))\ind_{[s\geq t']}ds
+\sum_{k=1}^n C(t_k,y'(t_k),\xi_k)\ind_{[t_{k}\leq
T]}+G(y(T))]\\
&& \quad \leq \integ{t'}{T}\psi(s,y'(s),\tau^1(s))\ind_{[s\geq
t']}ds +\sum_{k=1}^n C(t_k,y'(t_k),\xi_k)\ind_{[t_{k}\leq
T]}+G(y(T))+\epsilon^1.
\end{eqnarray*}

Then
\begin{eqnarray*}
&& q^n(t',x')-q^n(t,x)\leq
\integ{t'}{T}\psi(s,y'(s),\tau^1(s))\ind_{[s\geq t']}ds
+\sum_{k=1}^n C(t_k,y'(t_k),\xi_k)\ind_{[t_{k}\leq T]}\\&&
\qquad\qquad\qquad \qquad\qquad+G(y(T))-
\integ{t}{T}\psi(s,y(s),\tau^1(s))ds \\&&\qquad \qquad\qquad\qquad
\qquad-\sum_{k=1}^n C(t_k,y(t_k),\xi_k)\ind_{[t_{k}\leq T]}
-G(y(T))+2\epsilon^1.
\end{eqnarray*}
Next w.l.o.g we assume that $t'<t$. Then we deduce that:
\begin{equation}\label{con_sup}
\begin{array}{ll}q^n(t',x')-q^n(t,x)&\leq
\integ{t_0}{T}\displaystyle\{(\psi(s,y'(s),\tau^1(s))-\psi(s,y(s),\tau^1(s)))\ind_{[s\geq t]}\}ds\\
{}&\qquad +\integ{t_0}{T}\psi(s,y'(s),\tau^1(s))\ind_{[t'\leq s<
t]}ds \\{}&\qquad +\displaystyle\sum_{k=1}^n\{
C(t_k,y'(t_k),\xi_k)-C(t_k,y(t_k),\xi_k)\}\ind_{[t_{k}\leq
T]}+2\epsilon^1\\
{}&\leq
\integ{t_0}{T}\displaystyle\{|\psi(s,y'(s),\tau^1(s))-\psi(s,y(s),\tau^1(s))|\ind_{[s\geq
t]}\}\\{}&\qquad
+\integ{t_0}{T}|\psi(s,y'(s),\tau^1(s))|\ind_{[t'\leq s< t]}ds
\\
{}&\qquad +\displaystyle n\max\limits_{k=1}^n
|C(t_k,y'(t_k),\xi_k)-C(t_k,y(t_k),\xi_k)|+2\epsilon^1.
\end{array}
\end{equation}
Using the uniform continuity of $\psi$, $C$ in $y$ and property
(\ref{estimat2}), then the right-hand side of (\ref{con_sup}), the
first and the second term converges to 0 as $t'$ tends to $t$ and
$x'$ tends to $x$.\\
 Taking the limit as $(t',x')\rw (t,x)$ we obtain:
$$
\limsup_{(t',x')\rw (t,x)}q^n(t',x')\leq q^n(t,x) +2\epsilon^1.$$As
$\epsilon^1$ is arbitrary then sending $\epsilon^1\rightarrow 0,$ to
obtain:
$$
\limsup_{(t',x')\rw (t,x)}q^n(t',x')\leq q^n(t,x).$$ Therefore $q^n$
is
upper semi-continuous.\\
\textbf{Step 2. }  Now we show that $q^n$ is lower semi-continuous.\\
 Fix an arbitrary
$\epsilon_2>0$. Let $\varphi_2=((t_n)_{n\geq 1},(\xi_n)_{n\geq 1})$
belongs to $D_a$ such that
\begin{eqnarray*}
&&\sup\limits_{\tau\in \Omega}[\integ{t'}{T}\psi(s,y'(s),\tau(s))ds
+\sum_{k=1}^n C(t_k,y'(t_k),\xi_k)\ind_{[t_{k}\leq
T]}+G(y(T))] \\
&& \quad \leq \inf\limits_{\varphi_2\in  D_a}\sup\limits_{\tau\in
\Omega}[\integ{t}{T}\psi(s,y'(s),\tau(s))ds +\sum_{k=1}^n
C(t_k,y'(t_k),\xi_k)\ind_{[t_{k}\leq T]}G(y(T))]+\epsilon_2\\ &&
\quad =q^n(t',x')+\epsilon_2.
\end{eqnarray*}
Also, $$q^n(t,x)\leq \sup\limits_{\tau\in
\Omega}[\integ{t}{T}\psi(s,y(s),\tau(s))\ind_{[s\geq t]}ds
+\sum_{k=1}^n C(t_k,y(t_k),\xi_k)\ind_{[t_{k}\leq T]}+G(y(T))],$$
now, pick $\tau_2$ such that

\begin{eqnarray*}
&& \sup\limits_{\tau\in
\Omega}[\integ{t}{T}\psi(s,y(s),\tau(s))\ind_{[s\geq t]}ds
+\sum_{k=1}^n C(t_k,y(t_k),\xi_k)\ind_{[t_{k}\leq T]}+G(y(T))]\\
&& \qquad \leq \integ{t}{T}\psi(s,y(s),\tau_2(s))\ind_{[s\geq t]}ds
+\sum_{k=1}^n C(t_k,y(t_k),\xi_k)\ind_{[t_{k}\leq
T]}+G(y(T))+\epsilon_2.
\end{eqnarray*}

Then
\begin{eqnarray*}
&& q^n(t',x')-q^n(t,x)\geq \integ{t'}{T}\psi(s,y'(s),\tau_2(s))ds
+\sum_{k=1}^n C(t_k,y'(t_k),\xi_k)\ind_{[t_{k}\leq
T]}+G(y(T))\\&&\qquad\qquad \qquad\qquad-
\integ{t}{T}\psi(s,y(s),\tau_2(s))ds -\sum_{k=1}^n
C(t_k,y(t_k),\xi_k)\ind_{[t_{k}\leq T]}-G(y(T))\\&&\qquad\qquad
\qquad\qquad-2\epsilon_2.
\end{eqnarray*}
Next w.l.o.g we assume that $t'<t$. Then we deduce that:
\begin{equation}\label{con_sup1}
\begin{array}{ll}v(t',x')-v(t,x)&\geq
\integ{t_0}{T}\displaystyle\{(\psi(s,y'(s),\tau_2(s))-\psi(s,y(s),\tau_2(s)))\ind_{[s\geq
t]}\}ds\\{}&\qquad +\integ{t_0}{T}
\psi(s,y'(s),\tau_2(s))\ind_{[t'\leq s< t]}ds\\
{}&\qquad +\displaystyle\sum_{k=1}^n\{
C(t_k,y'(t_k),\xi_k)-C(t_k,y(t_k),\xi_k)\}\ind_{[t_{k}\leq
T]}-2\epsilon_2\\
{}&\geq
-\integ{t_0}{T}\displaystyle\{|\psi(s,y'(s),\tau_2(s))-\psi(s,y(s),\tau_2(s))|\ind_{[s\geq
t]}\}ds\\{}&\qquad
-\integ{t_0}{T}|\psi(s,y'(s),\tau_2(s))|\ind_{[t'\leq s< t]}ds\\
{}&\qquad -\displaystyle n\max\limits_{1\leq k\leq n}
|C(t_k,y'(t_k),\xi_k)-C(t_k,y(t_k),\xi_k)|-2\epsilon_2.
\end{array}
\end{equation}
Using the uniform continuity of $\psi$, $C$ in $y$ and property
(\ref{estimat2}). Then the right-hand side of (\ref{con_sup}) the
first and the second term converges to 0 as $t'\rightarrow t$ and
$x'\rightarrow x$.\\
Taking the limit as $(t',x')\rw (t,x)$ to obtain:
$$
\liminf_{(t',x')\rw (t,x)}q^n(t',x')\geq q^n(t,x) -2\epsilon_2.$$As
 $\epsilon_2$ is arbitrary then putting
$\epsilon_2\rightarrow 0$ to obtain:
$$
\liminf_{(t',x')\rw (t,x)}q^(t',x')\geq q^n(t,x).$$ Therefore $q^n$
is
lower semi-continuous. We then proved that $q^n$ is  continuous.\\
\textbf{Step 3. } Let us show that $v$ is continuous. we have:
\begin{equation}\label{ine-con}|v(t,x)-v(t',x')|\leq
|v(t,x)-q^n(t,x)|+|q^n(t,x)-q^n(t',x')|+|q^n(t',x')-v(t',x')|\end{equation}
For n large enough, then we put $ n \rightarrow +\infty $ and using
$\lim\limits_{n\rightarrow \infty}q^n(t,x)=v(t,x)$ and the
continuity of $q^n(t,x)$ in $t$ and x, we get that the right hand
side terms of (\ref{ine-con}) converge to 0 as
$(t^{'},x^{'})\rightarrow (t,x)$. Therefore $v(t',x') \rightarrow
v(t,x)$ as $(t^{'},x^{'})\rightarrow (t,x)$. So  $v$ is continuous.
$\Box$

\subsection{Terminal value}
Because of the possible jumps at the terminal time $T$, it is easy
to see that, in general, $v(t,x)$ does not tend to $G(x)$ as $t$
tends to $T$. Extend the set of jumps to include jumps of zero,
meaning no jump. Call this extended set $E_0$, extend trivially the
operator $N$ to a function independent from $t$, and let

\begin{equation}\label{Valeur-ter} G_1(x)=\inf\limits_{\xi \in
E_0}[G(x+g(T,x,\xi))+C(T,x,\xi)]=\min\{G(x),N[G](T,x)\}.
\end{equation}
We know that $G$ and $C$ are uniformly continuous in $x$ then
$G_1(x)$ is continuous. We claim
\begin{lem}

$$v(t,x)\rightarrow G_1(x)\quad \mbox{as} \quad t\rightarrow T.$$
\end{lem}
$Proof$: Fix $(t,x)$ and a strategy $\varphi$. As in the previous
proof, for each $\tau(·)$, gather all jumps of $\varphi(\tau)$ if
any, in jump $\xi_1$ at the time T. Then we have
$$|J(t,x,\varphi,\tau)-G(x+g(T,x,\xi_1))-C(T,x,\xi_1)|\leq C_x(T-t)$$
or
$$J(t,x,\varphi,\tau)=G(x+g(T,x,\xi_1))+C(T,x,\xi_1)+O(T-t).$$
The right hand side above only depends on $\xi_1$, not on $\tau(.)$
itself. It follows that
$$
\begin{array}{ll}\inf\limits_{\varphi}\sup\limits_{\tau}J(t,x,\varphi,\tau)&=\inf\limits_{\xi\in
E_0}[G(x+g(T,x,\xi))+C(T,x,\xi)]+O(T-t)\\ {}&=G_1(x)+O(T-t).
\end{array}$$
The result follows letting $t\rightarrow T.\Box$
\section{Viscosity characterization of the value function}

In this section we prove that the value function $v$ is a viscosity
solution of the Hamilton-Jacobi-Isaacs quasi-variational inequality,
that we replace by an equivalent QVI easier to investigate.\\

We now consider the following quasi-variational inequality (Isaacs
equation): \be
\begin{array}{l}\label{Issacs}
\displaystyle \max\left\{\min\limits_{\tau\in{\cal
K}}\left[-\frac{\partial v}{\partial t}- \frac{\partial v}{\partial
x}f(t,x,\tau)-\psi(t,x,\tau)\right],\right.\\\qquad\qquad \qquad
\left.v(t,x)-N[v](t,x)\right\}= 0,
\end{array}
\ee with the terminal condition: $v(t,x)=G_1(x),\, x\in \R^m$, where
$G_1$ is given by (\ref{Valeur-ter}).\\
Notice that it follows from hypothesis that the term in square
brackets in (\ref{Issacs}) above is continuous with respect to
$\tau$ so that the minimum in $\tau$ over the compact $\cal K$
exists.\\
Recall the notion of viscosity solution of QVI (\ref{Issacs}).
\begin{axiom} Let $v$ be a  continuous function defined on
$[0,T]\times \R^m$, $\R$-valued and such that $v(T,x)=G_1(x)$ for
any $x\in \R^m$. The $v$ is called:
\begin{itemize}
\item [$(i)$] A viscosity supersolution  of (\ref{Issacs})
if for any $(\overline{t},\overline{x})\in [t_0,T[\times \R^m$ and
any function $\varphi \in C^{1,2}([t_0,T[\times \R^m)$ such that
$\varphi(\overline{t},\overline{x})=v(\overline{t},\overline{x})$
and $(\overline{t},\overline{x})$ is a local maximum of $\varphi
-v$, we have: \be
\begin{array}{l}\label{Issacs1}
\displaystyle \max\left\{\min\limits_{\tau\in{\cal
K}}\left[-\frac{\partial \varphi}{\partial t}- \frac{\partial
\varphi}{\partial
x}f(\overline{t},\overline{x},\tau)-\psi(\overline{t},\overline{x},\tau)\right],\right.\\\qquad\qquad
\qquad
\left.v(\overline{t},\overline{x})-N[v](\overline{t},\overline{x})\right\}\geq
0.
\end{array}
\ee

\item [$(ii)$] A viscosity  subsolution of (\ref{Issacs})
if for any $(\overline{t},\overline{x})\in [t_0,T[\times \R^m$ and
any function $\varphi \in C^{1,2}([t_0,T[\times \R^m)$ such that
$\varphi(\overline{t},\overline{x})=v(\overline{t},\overline{x})$
and $(\overline{t},\overline{x})$ is a local minimum of $\varphi
-v$, we have: \be
\begin{array}{l}\label{Issacs2}
\displaystyle \max\left\{\min\limits_{\tau\in{\cal
K}}\left[-\frac{\partial \varphi}{\partial t}- \frac{\partial
\varphi}{\partial
x}f(\overline{t},\overline{x},\tau)-\psi(\overline{t},\overline{x},\tau)\right],\right.\\\qquad\qquad
\qquad
\left.v(\overline{t},\overline{x})-N[v](\overline{t},\overline{x})\right\}\leq0.
\end{array}
\ee
\item [$(iii)$] A viscosity solution if it is both a viscosity supersolution and
subsolution. $\Box$
\end{itemize}
\end{axiom}
\begin{theo}
The value function $v$ is the viscosity solution of the
quasi-variational inequality (\ref{Issacs}).
\end{theo}
$Proof$: The viscosity property follows from the dynamic programming
principle and is proved in \cite{[EBB]}.$\Box$\\

Now we give an equivalent of quasi-variational inequality
(\ref{Issacs}). In this section, we consider the new function
$\Gamma$ given by the classical change of variable $\Gamma(t,x) =
\exp(t)v(t, x)$, for any $t\in[t_0,T ]$ and $x\in \R^m$. Of course,
the function $\Gamma$ is bounded from below and continuous with
respect to its arguments.\\ A second property is given by the

\begin{pro}
$v$ is a viscosity solution of (\ref{Issacs}) if and only if
$\Gamma$ is a viscosity solution to the following quasi-variational
inequality in $[t_0,T [\times \R^m$, \be
\begin{array}{l}\label{Issacs3}
\displaystyle \max\left\{\min\limits_{\tau}\left[-\frac{\partial
\Gamma}{\partial t}+\Gamma(t,x)- \frac{\partial \Gamma}{\partial
x}f(t,x,\tau)-\exp(t)\psi(t,x,\tau)\right],\right.\\\qquad\qquad
\qquad \left.\Gamma(t,x)-M[\Gamma](t,x)\right\}= 0,
\end{array}
\ee where $M[\Gamma](t,x)=\inf\limits_{\xi\in
E}[\Gamma(t,x+g(t,x,\xi))+\exp(t)C(t,x,\xi)].$ The terminal
condition for $\Gamma$ is: $\Gamma(T,x)=\exp(T)G_1(x)$ in
$\R^m.$$\Box$
\end{pro}

\section{Uniqueness of the solution of quasi-variational inequality} We are going now to address the
question of uniqueness of the viscosity solution of
quasi-variational inequality (\ref{Issacs}). We have the following:

\beth \label{uni1}The solution in viscosity sense of
quasi-variational inequality (\ref{Issacs}) is unique in the space
of continuous functions on $[t_0,T]\times R^m$ which satisfy a
linear growth condition, i.e., in the space
$$\begin{array}{l}{\cal C}:=\{\varphi: [0,T]\times \R^m\rightarrow
\R, \mbox{ continuous and for any }\\\qquad \qquad\qquad(t,x), \,
\varphi(t,x)\leq C(1+|x|) \mbox{ for some constants } C \,\mbox{and
bounded from below} \}.\end{array}$$\eeth {\it Proof}. The proof is
divided in four steps. We will show by contradiction that if $u$ and
$w$ is a subsolution and a supersolution respectively for
(\ref{Issacs3}) then  $u\leq w$. Therefore if we have two solutions
of (\ref{Issacs3}) then they are obviously equal.\\

Now we fix $\lambda \in (0, 1)$, close to 0, and prove the
comparison result for $(1-\lambda)u$ and $w$. Actually for some
$R>0$ suppose there exists
$(\overline{t},\overline{x})\in(0,T)\times B_R $ $(B_R := \{x\in
\R^m; |x|<R\})$ such that:
\begin{equation}
\label{comp-equ1}
\max\limits_{t,x}((1-\lambda)u(t,x)-w(t,x))=(1-\lambda)u(\overline{t},\overline{x})-w(\overline{t},\overline{x})=\eta>0.
\end{equation}\textbf{Step 1. } Let us take $\beta >0$ and $\theta \in (0,1]$ small enough, so that the following holds:
\begin{equation}
\left\{
\begin{array}{l}
 - \beta
(T-\overline{t})^2-\lambda\exp(T)G_1(\overline{x})<\dfrac{\eta }{2} \\
2\theta \mid g(\overline{t},\overline{x},\xi)\mid ^4 <\lambda \exp(\overline{t})C(\overline{t},\overline{x},\xi), \quad \forall \xi \in \R^m .%\\
\end{array}%
\label{g} \right.
\end{equation}%. %, so that the
%following holds:
%\begin{equation}\label{beta} \beta
%=(\frac{-\lambda\exp(T)G_1(x)}{(T-\overline{t})^2}\vee 0).
%\end{equation}
 Then, for a small $\epsilon >0$, let us define:
\begin{equation}
\begin{array}{ll}
\label{phi}
\Phi_{\epsilon}(t,x,y)=(1-\lambda)u(t,x)-w(t,y)-\frac{1}{2\epsilon}|x-y|^{2}
-\theta(|x-\overline{x}|^{4}+|y-\overline{x}|^{4})\\
\qquad\qquad\qquad\qquad-\beta (t-\overline{t})^2.
\end{array}
\end{equation}
By the linear growth assumption on $u$ and $w$, there exists a
$(t_{0},x_{0},y_{0})\in [0,T]\times B_R \times B_R $, for $R$ large
enough, such that:
$$\Phi_{\epsilon}(t_{0},x_{0},y_{0})=\max\limits_{(t,x,y)}\Phi_{\epsilon}(t,x,y).$$
On the other hand, from $2 \Phi_{\epsilon}(t_{0},x_{0},y_{0})\geq
\Phi_{\epsilon}(t_{0},x_{0},x_{0})+\Phi_{\epsilon}(t_{0},y_{0},y_{0})$,
we have
\begin{equation}
\frac{1}{\epsilon}|x_{0} -y_{0}|^{2} \leq
(1-\lambda)(u(t_{0},x_{0})-u(\overline{t},\overline{x}))+(w(t_{0},y_{0})-w(\overline{t},\overline{x})),
\end{equation}
and consequently $\frac{1}{\epsilon}|x_{0} -y_{0}|^{2}$ is bounded,
and as $\epsilon\rightarrow 0$, $|x_{0} -y_{0}|\rightarrow 0$. Since
$u$ and $w$ are uniformly continuous on $[0,T]\times
\overline{B}_R$, then $\frac{1}{\epsilon}|x_{0}
-y_{0}|^{2}\rightarrow 0$ as
$\epsilon\rightarrow 0.$\\
Since
 $\Phi_{\epsilon}(t_{0},x_{0},y_{0})\geq
\Phi_{\epsilon}(\overline{t},\overline{x},\overline{x})$, we have
\begin{equation}
(1-\lambda)u(\overline{t},\overline{x})-w(\overline{t},\overline{x})\leq
\Phi_{\epsilon}(t_{0},x_{0},y_{0})\leq
(1-\lambda)u(t_{0},x_{0})+w(t_{0},y_{0}).
\end{equation}
it follow from the continuity of $u$ and $w$ that, up to a
subsequence,
\begin{equation}
\begin{array}{ccc}
 \label{subsequence}
 (t_0,x_0,y_0)\rightarrow (\overline{t},\overline{x},\overline{x})\\
 \theta(|x_0-\overline{x}|^{4}+|y_0-\overline{x}|^{4})\rightarrow 0\\
(1-\lambda)u(t_{0},x_{0})+w(t_{0},y_{0})\rightarrow
(1-\lambda)u(\overline{t},\overline{x})-w(\overline{t},\overline{x})
 \end{array}
 \end{equation}
\textbf{Step 2. } We now show that $t_{0} <T.$ Actually if $t_{0}
=T$ then,
$$
\Phi_{\epsilon}(\overline{t},\overline{x},\overline{x})\leq
\Phi_{\epsilon}(T,x_{0},y_{0}),$$ and,
$$
(1-\lambda)u(\overline{t},\overline{x})-w(\overline{t},\overline{x})\leq
(1-\lambda)\exp(T)G_1(x_{0}) -\exp(T)G_1(y_{0})- \beta
(T-\overline{t})^2,
$$
since $u(T,x_{0})=\exp(T)G_1(x_{0})$, $w(T,y_{0})=\exp(T)G_1(y_{0})$
and $G_1$ is uniformly continuous on $\overline{B}_R$. Then as
$\epsilon \rightarrow 0$ and (\ref{g}), we have,
$$
\begin{array}{ll}
\eta &\leq - \beta
(T-\overline{t})^2-\lambda\exp(T)G_1(\overline{x})\\
\eta & \leq \frac{\eta}{2},
\end{array}
$$
which yields a contradiction and we have $t_0 \in [0,T)$.\\
\textbf{Step 3. } We now claim that:
\begin{equation}
\label{visco-comp11}w(t_0,y_0)-\inf\limits_{\xi\in
E}[w(t_0,y_0+g(t_0,y_0,\xi))+\exp(t_0)C(t_0,y_0,\xi)] < 0.
\end{equation}
Indeed if
$$w(t_0,y_0)-\inf\limits_{\xi\in
E}[w(t_0,y_0+g(t_0,y_0,\xi))+\exp(t_0)C(t_0,y_0,\xi)]\geq 0,$$ then
from Proposition 2 we have:
$$w(t_0,y_0)-\inf\limits_{\xi\in
E}[w(t_0,y_0+g(t_0,y_0,\xi))+\exp(t_0)C(t_0,y_0,\xi)]= 0,$$

then there exist $\xi_1\in E$ and small $\epsilon_1>0$ such that:
$$w(t_0,y_0)-w(t_0,y_0+g(t_0,y_0,\xi_1))-\exp(t_0)C(t_0,y_0,\xi_1)\geq
-\epsilon_1.$$ From the subsolution property of $u(t_0,x_0)$, we
have
$$ u(t_0,x_0)-\inf\limits_{\xi\in
E}[u(t_0,x_0+g(t_0,x_0,\xi))+\exp(t_0)C(t_0,x_0,\xi)]\leq 0,
$$
then
$$u(t_0,x_0)-u(t_0,x_0+g(t_0,x_0,\xi_1))-\exp(t_0)C(t_0,x_0,\xi_1)\leq0.$$ It
follows that:
$$\begin{array}{ll}(1-\lambda)u(t_0,x_0)- w(t_0,y_0) -[(1-\lambda)u(t_0,x_0+g(t_0,x_0,\xi_1))
-w(t_0,y_0+g(t_0,y_0,\xi_1))] \\\leq
(1-\lambda)\exp(t_0)C(t_0,x_0,\xi_1)- \exp(t_0)C(t_0,y_0,\xi_1)
+\epsilon_1.\end{array}$$ Now taking into account of (\ref{phi}) to
obtain:
\begin{equation}
\begin{array}{l}
\Phi_{\epsilon}(t_{0},x_{0},y_{0})-\Phi_{\epsilon}(t_{0},x_{0}+g(t_0,x_0,\xi_1),y_{0}+g(t_0,y_0,\xi_1)) \\
\quad \leq -\frac{1}{2\epsilon}|x_0-y_0|^{2}
-\theta(|x_0-\overline{x}|^{4}+|y_0-\overline{x}|^{4})-\beta
(t_0-\overline{t})^2 \\
\quad +
\frac{1}{2\epsilon}|x_0+g(t_0,x_0,\xi_1)-y_0-g(t_0,y_0,\xi_1)|^{2}
 \\
\quad
+\theta(|x_0+g(t_0,x_0,\xi_1)-\overline{x}|^{4}+|y_0+g(t_0,y_0,\xi_1)-\overline{x}|^{4})+\beta
(t_0-\overline{t})^2 \\ \quad +
(1-\lambda)\exp(t_0)C(t_0,x_0,\xi_1)-
\exp(t_0)C(t_0,y_0,\xi_1)+\epsilon_1\\
\quad \leq
\frac{1}{2\epsilon}|x_0-y_0|^{2}+\theta|g(t_0,x_0,\xi_1)|^{4}+\theta|g(t_0,y_0,\xi_1)|^{4}\\
\quad  -\lambda \exp(t_0)C(t_0,x_0,\xi_1)+\exp(t_0)C(t_0,x_0,\xi_1)-
\exp(t_0)C(t_0,y_0,\xi_1)+\epsilon_1.%
\end{array}
\label{super-ine}
\end{equation}%
Here in the last inequality we have used $g$ is Lipschitz in $x$.
Now since  $C$ and $g$ are uniformly continuous on $[0,T]\times
\overline{B}_R$ then
\begin{equation}
\begin{array}{l}
\Phi_{\epsilon}(t_{0},x_{0},y_{0})-\Phi_{\epsilon}(t_{0},x_{0}+g(t_0,x_0,\xi_1),y_{0}+g(t_0,y_0,\xi_1))\leq\\
\qquad\qquad\qquad\qquad
 -\lambda \exp(\overline{t})C(\overline{t},\overline{x},\xi_1)+2\theta \mid g(\overline{t},\overline{x},\xi_1)\mid ^4+\epsilon(2+2\theta+\lambda)+\epsilon_1  .%
\end{array}
\label{super-ine1}
\end{equation}
Since (\ref{g}), sending $\epsilon\rightarrow0$ and sending
$\epsilon_1\rightarrow0$, we get
$$\Phi_{\epsilon}(t_{0},x_{0},y_{0})<\Phi_{\epsilon}(t_{0},x_{0}+g(t_0,x_0,\xi_1),y_{0}+g(t_0,y_0,\xi_1))$$
This is contradiction to the fact that $(t_0,x_0,y_0)$ is the
supremum point of $\Phi_{\epsilon}$. Then the claim
(\ref{visco-comp11}) holds.\\

\textbf{Step 4. } To complete the proof it remains to show
contradiction. Let us denote
\begin{equation}
\varphi_{\epsilon}(t,x,y)=\frac{1}{2\epsilon}|x-y|^{2}
+\theta(|x-\overline{x}|^{4}+|y-\overline{x}|^{4})+\beta
(t-\overline{t})^2.
\end{equation}
Then we have: \be \left\{
\begin{array}{lllll}\label{derive}
D_{t}\varphi_{\epsilon}(t,x,y)=2\beta(t-\overline{t}),\\
D_{x}\varphi_{\epsilon}(t,x,y)= \frac{1}{\epsilon}(x-y) +4\theta
(x-\overline{x})|x-\overline{x}|^{2}, \\
D_{y}\varphi_{\epsilon}(t,x,y)= -\frac{1}{\epsilon}(x-y) +
4\theta(y-\overline{x})|y-\overline{x}|^{2}.
\end{array}
\right. \ee Let $c,d\in \R$ such that
\begin{equation}\label{uni-vis}c+d=2\beta(t_0-\overline{t}).\end{equation}
 Taking now into account
(\ref{visco-comp11}), (\ref{uni-vis}), and the definition of
viscosity solution, we get:
\begin{equation}\begin{array}{lll}\label{vis_sub1}\displaystyle\min\limits_{\tau}[-c+(1-\lambda)u(t_0,x_0)
-\langle\frac{1}{\epsilon}(x_0-y_0) +4\theta
(x_0-\overline{x})|x_0-\overline{x}|^{2},\\\qquad\qquad\qquad\qquad
f(t_0,x_0,\tau)\rangle-(1-\lambda)\exp(t_0)\psi(t_0,x_0,\tau)]\leq0
\end{array}\end{equation}

 and

\begin{equation}\begin{array}{l}\label{vis_sub11}\displaystyle\min\limits_{\tau}\left[d+w(t_0,y_0)-\langle
\frac{1}{\epsilon}(x_0-y_0) -4\theta
(y_0-\overline{x})|y_0-\overline{x}|^{2},
f(t_0,y_0,\tau)\rangle-\exp(t_0)\psi(t_0,y_0,\tau)\right]\\\geq0,
\end{array}\end{equation}

which implies that:
\begin{equation}
\begin{array}{llllll}
\label{viscder} &-c-d+(1-\lambda)u(t_0,x_0)-w(t_0,y_0)\\& \leq
\displaystyle\min\limits_{\tau}\left[-\langle
\frac{1}{\epsilon}(x_0-y_0) -4\theta
(y_0-\overline{x})|y_0-\overline{x}|^{2},
f(t_0,y_0,\tau)\rangle -\exp(t_0)\psi(t_0,y_0,\tau)\right]\\
&
-\displaystyle\min\limits_{\tau}\left[-\langle\frac{1}{\epsilon}(x_0-y_0)
+4\theta (x_0-\overline{x})|x_0-\overline{x}|^{2},
f(t_0,x_0,\tau)\rangle-(1-\lambda)\exp(t_0)\psi(t_0,x_0,\tau)\right],\\&
\leq
\displaystyle\sup\limits_{\tau}\left[\langle\frac{1}{\epsilon}(x_0-y_0)
+4\theta (x_0-\overline{x})|x_0-\overline{x}|^{2},
f(t_0,x_0,\tau)\rangle+(1-\lambda)\exp(t_0)\psi(t_0,x_0,\tau)\right]\\
&-\displaystyle\sup\limits_{\tau}\left[\langle
\frac{1}{\epsilon}(x_0-y_0) -4\theta
(y_0-\overline{x})|y_0-\overline{x}|^{2}, f(t_0,y_0,\tau)\rangle
+\exp(t_0)\psi(t_0,y_0,\tau)\right],\\& \leq
\displaystyle\sup\limits_{\tau}[\langle\frac{1}{\epsilon}(x_0-y_0) ,
f(t_0,x_0,\tau)-f(t_0,y_0,\tau)\rangle\\&+ \langle 4\theta
(x_0-\overline{x})|x_0-\overline{x}|^{2}, f(t_0,x_0,\tau)\rangle +
\langle4\theta (y_0-\overline{x})|y_0-\overline{x}|^{2},
f(t_0,y_0,\tau)\rangle\\
&

+(1-\lambda)\exp(t_0)\psi(t_0,x_0,\tau)-\exp(t_0)\psi(t_0,y_0,\tau)].
\end{array}
\end{equation}

Now, from (\ref{regbs1}), we get:
$$
\langle\frac{1}{\epsilon}(x_0-y_0),f(t_0,x_0,\tau)-f(t_0,y_0,\tau)\rangle
\leq \frac{C}{\epsilon}|x_0 - y_0|^{2}.$$ Next $$ \langle 4\theta
(x_0-\overline{x})|x_0-\overline{x}|^{2}, f(t_0,x_0,\tau)\rangle\leq
4C\theta |x_0||x_0-\overline{x}|^{3},$$and finally,
$$
\langle 4\theta (y_0-\overline{x})|y_0-\overline{x}|^{2},
f(t_0,y_0,\tau)\rangle\leq 4C\theta
|y_0||y_0-\overline{x}|^{3}.$$Taking in to account
$$c+d=2\beta(t_0-\overline{t}).$$ So that by plugging into (\ref{viscder}) and note that $\lambda
>0$ we obtain:
\begin{equation}
\begin{array}{llll}
\label{viscder}
&-2\beta(t_0-\overline{t})+(1-\lambda)u(t_0,x_0)-w(t_0,y_0)\\&
\qquad\qquad\qquad\qquad\leq \displaystyle\frac{C}{\epsilon}|x_0 -
y_0|^{2}\\&\qquad\qquad\qquad\qquad\qquad+ 4C\theta
|x_0||x_0-\overline{x}|^{3} + 4C\theta
|y_0||y_0-\overline{x}|^{3}\\
&\qquad\qquad\qquad\qquad\qquad

+\sup\limits_{\tau}(1-\lambda)\exp(t_0)\psi(t_0,x_0,\tau)-\exp(t_0)\psi(t_0,y_0,\tau)].
\end{array}
\end{equation}By sending $\epsilon \rightarrow0$, $\theta
\rightarrow0$ and taking into account of the continuity of $\psi$,
we obtain $\eta \leq 0$ which is a contradiction. Now sending
$\lambda\rightarrow0$, we get the required comparison between $u$
and $w$. The proof of
Theorem \ref{uni1} is now complete. $\Box$\\

%As an example of a use of this result in mathematical finance, one
%may consider the option pricing problem of references \cite{[BE],
%[BET]}.  If the piecewise linear transaction costs are replaced by a
%more realistic piecewise affine cost, i.e. a fixed cost is charged
%for any transaction in addition to a variable part, then the problem
%at hand is exactly that considered here. %This was actually the
%motivation for the present analysis.
\medskip

%\no {\bf Acknowledgement}: The author thanks gratefully Prof. A.
%Popier for the fructuous discussions during the preparation of this
%paper. I also would like to thank the referees for their careful
%reading and for their helpful comments and suggestions that led to
%considerable improvements in the paper.$\Box$

%As a by-product we have the following corollary: \bcor Let
%$(v^1,...,v^m)$ be a viscosity solution of (\ref{sysvi1}) which
%satisfies a polynomial growth condition then for $i=1,...,m$ and
%$(t,x)\in [0,T]\times \R^k$, $$ v^i(t,x)= \sup_{(\delta,\xi)\in
%{\cal D}^i_t}E[\integ{t}{T}\psi_{u_s}(s,X^{tx}_s)ds -\sum_{n\geq 1}
%g_{u_{\tau_{n-1}}u_{\tau_n}}(\tau_{n},X^{tx}_{\tau_{n}})
%\ind_{[\tau_{n}<T]}].
%$$
%\section {Numerical results} We consider now some numerical examples
%of
%the optimal switching problem (\ref{sysvi1}).\\
%{\bf Example1}: In this example we consider an optimal switching
%problem with two
%modes, where\\
%$T=1$, $b=x$, $\sigma=\sqrt{2}x$, $g_{12}(t,x)=0.5|x|+1$,
%$g_{21}(t,x)=0.1|x|+0.5t+1$, $\psi_1(t,x)=x+0.75t+1$,
%$\psi_2(t,x)=0.1x+t-2$.
%
%
% \ecor
%
%\no {\bf Acknowledgement}: the authors thank gratefully Prof.
%J.Zhang for the fructuous discussions during the preparation of this
%paper.$\Box$


\begin{thebibliography}{99}
\small
\renewcommand{\baselinestretch}{0.3}

\bibitem{[B]} Barles, G. (1985): Deterministic impulse control problems. {\it SIAM J.
Control Optim.} (23), pp. 419-432.

\bibitem{[BEJ]} Barron, E.N., Evans, L.C. and  Jensen, R. (1984): Viscosity solutions of Isaaes' equations and differential games with
Lipschitz controls. {\it J Differential Equations.} (53), pp.
213-233.

\bibitem{[BL]} Bensoussan, A. and Lions, J.L. (1984): Impulse Control and Quasi-Variational
Inequalities. {\it Bordes, Paris.}

\bibitem{[BE]} Bernhard, P. (2005): A robust control approach to option pricing including
transaction costs. {\it Annals of the ISDG.}(7), pp. 391-416.
Birkhaüser, Basel.

\bibitem{[BET]}Bernhard, P., El Farouq, N. and Thiery, S. (2006): An impulsive differential
game arising in finance with interesting singularities. {\it Annals
of the ISDG.}  (8), pp. 335-363. Birkhaüser, Basel.


\bibitem{[CE]} Capuzzo-Dolcetta, I. and  Evans,  L.C. (1984): Optimal switching for ordinary
differential equations, SIAM J. Control Optim. (22), pp 143-161.

\bibitem{[CIL]} Crandall, M., Ishii, H. and P.L. Lions (1992): User's guide to viscosity solutions of
second order partial differential equations, Bull. Amer. Math. Soc.,
27, 1-67.

\bibitem{[DS]} Dharmatti, S. and Shaiju, A.J. (2007): Infinite
dimensional differential games with hybrid controls. {\it Proc.
Indian Acad. Sci. Math.} (117), pp. 233-257.

\bibitem{[DR]}Dharmatti, S. and Ramaswamy, M. (2006): Zero-sum differential games involving
hybrid controls. {\it J. Optim. Theory Appl.} (128), pp. 75-102.


\bibitem{[E]} El Asri, B. (2010): Optimal Multi-Modes Switching Problem in Infinite Horizon. {\it Stochastics and
Dynamics.} (10), pp. 231-261.

\bibitem{[E1]} El Asri, B. (2013): Deterministic minimax impulse control in finite horizon: the viscosity solution approach.
 {\it ESAIM: Control, Optimisation and Calculus of Variations.} (19), pp. 63-77.

\bibitem{[E2]} El Asri, B. (2013): Stochastic optimal multi-modes switching with a viscosity solution approach. {\it Stochastic Processes and their Applications.}
 (123), pp. 579-602.



\bibitem{[EBB]}  El Farouq, N., Barles, G. and Bernhard, P. (2010): Deterministic Minimax Impulse Control
{\it Appl.Math. Optim.} DOI 10.1007/s00245-009-9090-0.

\bibitem{[ES]} Evans, L.C., Souganidis, P.E. (1984): Differential games and representation
formulas for the solution of Hamilton-Jacobi-Isaacs equations. {\it
Indiana Univ. J. Math.} (33), pp. 773-797.

\bibitem{[F]}Fleming, W.H. (1964): The convergence problem for differential games, 2.
{\it Ann. Math. Study} (52), pp. 195-210.

 \bibitem{[L]}Lions, P.L. (1982): Generalized
Solutions of Hamilton-Jacobi Equations. {\it Pitman, London.}

\bibitem{[LS]} Lions, P.L. and Souganidis, P.E. (1985): Differential games, optimal control
and directional derivatives of viscosity solutions of Bellman's and
Isaacs' equations. {\it SIAM J. Control Optim.} (23), pp. 566-583
\bibitem{[SD]} Shaiju, A.J. and Dharmatti, S. (2005): Differential games with
continuous, switching and impulse controls. {\it Nonlinear Anal.}
(63), pp. 23-41.
\bibitem{[S]}Souganidis, P.E. (1985): Max-min representations and
product formulas for the viscosity solutions of Hamilton-Jacobi
equations with applications to differential games. {\it Nonlinear
Anal. Theory Methods Appl.} (9), pp. 217-57.

\bibitem{[Y11]}Yong, J.M. (1989): Systems governed by ordinary differential equations
with continuous, switching and impulse controls. {\it Appl Math
Opti.} (20), pp. 223-235..

\bibitem{[Y1]}Yong, J.M. (1989): Optimal switching and impulse controls for
distributed parameter systems. {\it Systems Sci Math Sci} (2), pp.
137-160. 23.
\bibitem{[Y2]}Yong, J.M. (1990): Differential games with switching strategies. {\it J Math
Anal Appl.} (145), pp. 455-469.
\bibitem{[Y3]}Yong, J.M. (1990): A zero-sum differential game in a
finite duration with switching strategies. {\it SIAM J Control
Optim.} (28), pp. 1234-1250.

\bibitem{[Y]}Yong, J.M. (1994): Zero-sum
differential games involving impulse controls. {\it Appl.Math.
Optim.} (29), pp. 243-261.



%\bibitem{[BS]} Brennan, M. J. and Schwartz, E. S. (1985): Evaluating
%natural resource investments. {\it J.Business} 58,  pp. 135-137.
%
%\bibitem{[CK]} Cvitanic, J. and Karatzas, I (1996): Backward SDEs with
%reflection and Dynkin games. {\it Annals of Probability} 24 (4), pp.
%2024-2056.
%
%\bibitem{[CL]} Carmona, R. and Ludkovski, M. (2005): Optimal
%Switching with Applications to Energy Tolling Agreements. {\it
%Preprint}.
%
%
%\bibitem{[DX]} Deng, S. J. and Xia, Z. (2005): Pricing and hedging
%electric supply contracts: a case with tolling agreements. {\it
%Preprint}.
%
%\bibitem{[D]} Dixit, A. (1989): Entry and exit decisions under
%uncertainty. {\it J. Political Economy} 97,  pp. 620-638.
%
%\bibitem{[DP]} Dixit, A. and  Pindyck, R. S. (1994): Investment under
%uncertainty. {\it Princeton Univ. Press.}
%
%\bibitem{[DH]} Djehiche, B. and Hamad\`ene, S (2007): On a finite horizon
%Starting and Stopping Problem with Default risk. {\it Preprint,
%Universit\'e du Maine, F.}
%
%\bibitem{[DHP]} Djehiche, B., Hamad\`ene, S. and Popier, A (2007): A finite horizon optimal multiple switching problem.
%{\it Preprint, Universit\'e du Maine, F.}
%
%
%\bibitem{[DZ]} Duckworth, K. and Zervos, M. (2000): A problem of
%stocahstic impulse control with discretionary stopping. In
%{Proceedings of the 39th IEEE Conference on Decision and Control},
%IEEE Control Systems Society, Piscataway, NJ, pp. 222-227.
%
%\bibitem{[DZ2]} Duckworth, K. and Zervos, M. (2001): A model for
%investment decisions with switching costs. {\it Annals of Applied
%probability} 11 (1), pp. 239-260.
%
%\bibitem{Elka} El Karoui, N. (1980): Les aspects probabilistes du
%contr\^ ole stochastique. {\it Ecole d'\'et\'e de probabilit\'es de
%Saint-Flour, Lect. Notes in Math. No 876, Springer Verlag.}
%
%
%\bibitem{[EKal]} El-Karoui, N., Kapoudjian, C., Pardoux, E., Peng, S.
%and Quenez, M. C. (1997): Reflected solutions of backward SDEs and
%related obstacle problems for PDEs. {\it Annals of Probability} 25
%(2), pp. 702-737.
%
%\bibitem{ham} Hamad\`ene, S. (2002): Reflected BSDEs with discontinuous
%barriers. {\it Stochastics and Stochastic Reports} 74 (3-4), pp.
%571-596.
%
%
%\bibitem{[HJ]} Hamad\`ene, S. and  Jeanblanc, M (2007):
%On the Starting and Stopping Problem: Application in reversible
%investments, {\it Math. of Operation Research, vol.32, No.1,
%pp.182-192}.
%
%\bibitem{hib} Hamad\`ene, S. and Hdhiri, I. (2006): On the starting and
%stopping problem in the model with jumps. {\it Preprint ,
%Universit\'e du Maine, Le Mans, F.}
%
%\bibitem{ht} Hu, Y., Tang, S. (2007):
%Multi-dimensional BSDE with Oblique Reflection and Optimal
%Switching. {\it Preprint Universit\'e de Rennes 1, France}
%
%
%\bibitem{[LP]} Ly Vath, V. and Pham, H. (2007): Explicit solution to an optimal switching problem in the two-regime case.
%SIAM Journal on Control and Optimization, pp. 395-426.



\end{thebibliography}
\end{document}